\documentstyle[12pt]{article}

\begin{document}
\newcommand{\ol }{\overline}
\newcommand{\ul }{\underline }
\newcommand{\ra }{\rightarrow }
\newcommand{\lra }{\longrightarrow }
\newcommand{\ga }{\gamma }
\newcommand{\st }{\stackrel }
\newcommand{\scr }{\scriptsize }
\title{{\bf The Baer Invariant of Semidirect and Verbal Wreath Products of Groups}
\footnote{This research was in part supported by a grant from the
Research Council of Ferdowsi University of Mashhad (No.
241-1316).}}
\author{Behrooz Mashayekhy
\\ Department of Mathematics, Ferdowsi University
of Mashhad,\\ P.O.Box 1159-91775, Mashhad, Iran\\ E-mail:
mashaf@math.um.ac.ir}
\date{ }
\maketitle
\begin{abstract}
 W. Haebich (1977, Journal of Algebra {\bf 44}, 420-433) presented
some formulas for the Schur multiplier of a semidirect product
and also a verbal wreath product of two groups. The author (1997,
Indag. Math., (N.S.), {\bf 8}({\bf 4}), 529-535) generalized a
theorem of W. Haebich to the Baer invariant of a semidirect
product of two groups with respect to the variety of nilpotent
groups of class at most $c\geq 1,\ {\cal N}_c$. In this paper,
first, it is shown that ${\cal V}M(B)$ and ${\cal V}M(A)$ are
direct factors of ${\cal V}M(G)$, where $G=B\rhd\!\!\!<A$ is the
semidirect product of a normal subgroup $A$ and a subgroup $B$ and
$\cal V$ is an arbitrary variety. Second, it is proved that
${\cal N}_cM(B\rhd\!\!\!<A)$ has some homomorphic images of
Haebich's type. Also some formulas of Haebich's type is given for
${\cal N}_cM(B\rhd\!\!\!<A)$, when $B$ and $A$ are cyclic groups.
 Third, we will present a formula for the Baer invariant of a $\cal V$-verbal
wreath product of two groups with respect to the variety of
nilpotent groups of class at most $c\geq 1$, where $\cal V$ is an
arbitrary variety. Moreover, it is tried to improve this formula,
when $G=A{\it Wr_V}B$ and $B$ is cyclic.
 Finally, a structure for the Baer invariant of a free wreath product with
respect to ${\cal N}_c$ will be presented, specially for the free
wreath product $A{\it Wr_*}B$ where $B$ is a cyclic group.

 {\it A.M.S. Classification 2000} : 20E10, 20E22, 20F12

\end{abstract}

\begin{center}
{\bf 1. INTRODUCTION AND MOTIVATION}
\end{center}

 As stated in my previous paper (2001, Journal of Algebra {\bf 235}, 15-26),
there are some results on the {\it Schur multiplier} and the {\it
Baer invariant} of the direct product of groups by I. Schur [16],
J. Wiegold [18], M.R.R. Moghaddam [13], G.Ellis [4], and the
author [12].

 It is known that a semidirect product of two groups is a generalization of a
direct product of groups. Also, we know that a verbal wreath
product is a kind of semidirect product. Therefore, it is
interesting to find some formulas for the Schur multiplier or the
Baer invariant of a semidirect product and a verbal wreath
product.

 In 1972, K.I. Tahara [17, (Theorem 2.2.5)], using cohomological methods, gave
a structure for a semidirect product of two groups. W. Haebich
[6] in 1977, presented some formulas for the Schur multiplier of
a semidirect product and a $\cal V$-verbal wreath product of two
groups, where $\cal V$ is an arbitrary variety. His method was
based on presentations of groups.

 In 1972, N. Blackburn [2] gave an explicit formula for the Schur multiplier of
a standard wreath product of two groups. Also E.W. Read [15] in
1976, found an explicit formula for the Schur multiplier of a
wreath product (we know that a standard wreath product is a
special kind of a wreath product). It should be mentioned that
their formulas were more explicit than Haebich's one. But
Haebich's results were more general than Blackburn and Read's
formulas.

 In 1997, the author [10] generalized a result of W. Haebich [6, (Theorem 1.7)]
to the Baer invariant of a semidirect product with respect to the
variety ${\cal N}_c$.

 Now, in this paper, we concentrate on the Haebich's method in order
to find some structures for the Baer invariant of a semidirect
product and a verbal wreath product of two groups with respect to
the variety of nilpotent groups of class at most $c,\ {\cal N}_c$.

 More precisely, first, using the functorial property of the Baer invariant, we
show that ${\cal V}M(B)$ and ${\cal V}M(A)$ are direct factors of
${\cal V}M(G)$, where $G=B\rhd\!\!\!<A$ is the semidirect product
of $A$ by $B$ and $\cal V$ is an arbitrary variety of groups
(Theorem 3.2). Also, in Section 3, we find some homomorphic images
and in finite case subgroups with structures similar to Haebich's
type [6,(Theorems 2.2 and 2.3)] for the Baer invariant of a
semidirect product $B\rhd\!\!\!<A$ with respect to the variety
${\cal N}_c$ (Theorem 3.6). Moreover, if $B$ and $A$ are cyclic
groups, then we show that the above structures are isomorphic to
${\cal N}_cM(B\rhd\!\!\!<A)$ (Theorem 3.7 and Corollary 3.11).

 In Section 4, we concentrate on a verbal wreath product and first, we find
a structure for ${\cal N}_cM(A{\it Wr_V}B)$ similar to Haebich's
type (Theorem 4.2). Second, we improve this structure for a
homomorphic image and in finite case for a subgroup of ${\cal
N}_cM(A{\it Wr_V}B)$ (Theorem 4.5). Moreover, if $B$ is cyclic,
then we prove that the last structure is isomorphic to ${\cal
N}_cM(A{\it Wr_V}B)$ (Theorem 4.6).

 Finally, in Section 5, using the previous results, we try to find a
structure for the Baer invariant of a free wreath product with
respect to the variety ${\cal N}_c$. Specially, we present a
formula for ${\cal N}_cM(A{\it Wr_*}B)$, when $B$ is cyclic and
${\cal N}_cM(A{\it Wr_*}B)$ is finite (Theorem 5.5), and also a
formula when $A$ and $B$ are both cyclic groups
(Theorem 5.6). \\

\begin{center}
{\bf 2. NOTATION AND PRELIMINARIES}
\end{center}

 We assume that the reader is familiar with the notions of variety of groups,
verbal and marginal subgroups (see [14]). Let $G$ be a group with
a free presentation $G\cong F/R$. Then the Baer invariant of $G$
with respect to the variety $\cal V$, denoted by ${\cal V}M(G)$
is defined to be
$$ {\cal V}M(G)=\frac {R\cap V(F)}{[RV^*F]}\ , $$
where $V(F)$ is the verbal subgroup of $F$ and
$$[RV^*G]=<v(g_1,\ldots ,g_{i-1},g_ia,g_{i+1},\ldots ,g_n)v(g_1,\ldots
,g_i,\ldots g_n)^{-1}\ |\ a\in R, $$ $$1\leq i\leq n,v\in V,g_i\in G,n\in
{\bf N}>\ ,$$
 which is independent of the choice of the presentation of $G$ and it is always
an abelian group (see [1,5,9] for further properties of the Baer
invariant).

 In special case, if $\cal V$ is the variety of abelian groups, $\cal A$, then
the Baer invariant of a group $G$ will be
$$ \frac {R\cap F'}{[R,F]}\ , $$
the {\it Schur multiplier} of $G$, where, in finite case, by Hopf's formula [7]
is isomorphic to the {\it second cohomology group} of $G$.

 Also, if $\cal V$ is the variety of nilpotent groups of class at most $c\geq
1$, ${\cal N}_c$, then the Baer invariant of a group $G$ will be
$$ {\cal N}_cM(G)=\frac {R\cap \ga_{c+1}(F)}{[R,\ _cF]}\ ,$$
where $[R,\ _cF]=[R,\underbrace{F,F,\ldots ,F}_{c-times}]$. The above notion is
also called the {\it $c$-nilpotent multiplier} of $G$ (see [3]).

 In order to deal with the Baer invariant of a semidirect product we need a
presentation for the semidirect product which is given as follows.

{\bf L}{\scr\bf EMMA} {\bf 2.1} [6]\ \  Suppose $G$ is a semidirect product (or a
splitting extension) of $A$ by $B$ under $\theta :B\lra Aut(A)$ and
$$ 1\lra R_1\lra F_1\st{\nu_1}{\lra}A\lra 1\ \ \ ,\ \ \ 1\lra R_2\lra
F_2\st{\nu_2}{\lra}B\lra 1 $$
are free presentations for $A$ and $B$, respectively. Then
$$ 1\lra R\lra F\lra G\lra 1 $$
is a free presentation for $G$ , where

$(i)$ $F=F_1*F_2$ , the free product of $F_1$ and $F_2$ ;

$(ii)\ R=R_1^FR_2^FS\ ;$

$(iii)\ S=<f_1^{-1}\ol{f_1}[f_2,f_1] | f_1,\ol{f_1}\in F_1;f_2\in
F_2;\nu_1(\ol {f_1})=\theta(\nu_2(f_2))(\nu_1f_1)>^F.$\\

 W. Haebich [6] in 1977, using the above lemma prove the following theorem about
the Schur multiplier of a semidirect product.

{\bf T}{\scr\bf HEOREM} {\bf 2.2} [6]\ \ Let $G$ be a semidirect
product of $A$ by $B$ under $\theta $ and $F_1/R_1$ and $F_2/R_2$
be presentations for $A$ and $B$, respectively. Then, by the
above notation, the following isomorphism holds:
$$ M(G)\cong M(B)\oplus \frac {S\cap F'}{[R_2,F_1][S,F]}\ .$$

 A generalization of the above theorem was presented by the author
in 1997 as follows.

{\bf T}{\scr\bf HEOREM} {\bf 2.3} [10]\ \ Let $G$ be a semidirect product of $A$
by $B$
under $\theta :B\lra Aut(A)$ and ${\cal N}_c$ be the variety of nilpotent groups
of class at most $c$ $(c\geq 1)$. Then
$$ {\cal N}_cM(G)\cong {\cal N}_cM(B)\oplus \frac {S\cap \ga_{c+1}(F)}{\prod
[R_2,F_1,F_2]_c[S,\ _cF]}\ \ ,$$
where $$ \prod [R_2,F_1,F_2]_c=<[r_2,f_1.\ldots ,f_c]\ |\ f_i\in F_1\cup
F_2,1\leq i\leq c,\exists k\ f_k\in F_1>^F\ .$$
In particular, ${\cal N}_cM(B)$ can be regarded as a direct factor of ${\cal
N}_cM(G)$ .\\

 We also assume that the reader is familiar with notions of a verbal product and
the cartesian subgroup of a free product ( see [11,14]). We need the following
properties of the cartesian subgroup.

{\bf L}{\scr\bf EMMA} {\bf 2.4} [14]\ \ Let $\{A_i|i\in I\}$ be a family of
groups and $\prod_{i\in I}^{*}A_i$ be the free product of
$A_i$'s, and let $[A_i]^*$ be the cartesian subgroup of the above free product.
Then

$(i)$ The cartesian subgroup avoids the constituents i.e $[A_i]^*\cap A_j=1$ ,
for all $j\in I$.

$(ii)$ If $a\in \prod_{i\in I}^{*}A_i$, then $a=a_{i_1}a_{i_2}\ldots a_{i_m}c$ ,
where $a_{i_j}\neq 1$ for all $j$,\\ $i_1<i_2<\ldots <i_m$ and $c\in [A_i]^*$.
The elements $a_{i_j}$ and $c$ are uniquely determined by $a$ and chosen order
of $I$.

{\bf L}{\scr\bf EMMA} {\bf 2.5} [14]\ \ Let $\cal V$ be a variety of groups
defined by the set of laws $V$. Then

 $(i)$ If $V(A)$ is the verbal subgroup of $A=\prod_{i}^{*}A_i$ , then $V(A)\cap
A_i=V(A_i)$ for all $i\in I$.

 $(ii)$ If $A=\prod_{i}^{*}A_i$, then
$$ V(A)=(\prod_{i}^{}V(A_i))(V(A)\cap [A_i]^*)\ . $$

 Now, in the following, you can find the definition of a {\it verbal wreath
product}.

{\bf D}{\scr\bf EFINITION} {\bf 2.6}\ \ Given arbitrary groups $A$ and $B$, let $A_b$ be an isomorphic copy of $A$ for
each $b\in B$ and denote by $a_b$ the element of $A_b$ mapped to $a\in A$.
Consider the ${\cal V}$-verbal product of $A_b$'s, $G_V= \prod_{b\in
B}^{V}A_b$ corresponding to the variety $\cal V$ with a set of words $V$,
which is isomorphic to the quotient $C/C_V$, where $C=\prod_{b\in B}^{*}A_b$
and $C_V=V(C)\cap [A_b]^*$. The map $a_b\mapsto a_{bb'}$ for all $a\in A\ ,\
b\in B$ and fixed $b'\in B$ induces automorphisms $\theta _*b'$ and $\theta_Vb'$
of $C$ and $G_V$, respectively, i.e.
\newpage
$$ \theta_*\ :\ B\lra Aut(C)\ \ \ \ \ \ ,\ \ \ \ \ \ \ \ \theta_V\ :\ B\lra
Aut(G_V)$$
$$\ \ \ \ \ \ \ \ \ \ \ \ \ \ b'\mapsto \theta_*b': C\ra C\ \ \ \ \ \ \ \ \ \ \
\ \ \ \ \ \ \ \ b'\mapsto \theta_Vb': G_V\ra G_V $$
$$\ \ \ \ \ \ \ \ \ \ \ \ \ \ \ \ \ \ \ \ \ \ \ \ \ \ \ \  a_b\mapsto a_{bb'}\ \ \ \ \ \ \ \ \ \ \ \ \ \
\ \ \ \ \ \ \ \ \ \ \ \ \ \ \ \ \ {\ol a}_b\mapsto {\ol a}_{bb'}$$
 It can be proved that $\theta_*:B\ra Aut(C)$ and $\theta_V:B\ra Aut(G_V)$ are
monomorphisms.

 Now, the $\cal V$-{\it verbal wreath product} of $A$ by $B$ is the semidirect
product of $G_V$ by $B$ under $\theta_V$, denoted by $A{\it
Wr_V}B$. The two special cases of the $\cal V$-verbal wreath
product are important, {\it the free wreath product} and {\it the
standard wreath product}.

 $(i)$ If $\cal V$ is the variety of all groups, then the $\cal V$-verbal wreath
product is the free wreath product which is the semidirect
product of $C$ by $B$ under $\theta_*$, denoted by $A{\it Wr_*}B$.

 $(ii)$ If ${\cal V}={\cal A}$ is the variety of abelian groups then the $\cal
V$-verbal wreath product is the standard wreath product which is the semidirect
product of the direct product $\prod_{b\in B}^{\times }A_b$ by $B$ under
$$ \theta : B\lra Aut(\prod_{b\in B}^{}\!^{\times }A_b)$$
$$\ \ \ \ \ \ \ \ \ \ \ \ \ \ b'\mapsto \theta b':\prod_{b}^{}\!^{\times}A_b\lra
\prod_{b}^{}\!^{\times}A_b$$ $$\ \ \ \ \ \ \ \ \ \ \ \ \ \ \ \ \
\ \ \ \ \ \ \ \ \{a_b\}_b\mapsto \{a_{bb'}\}_b$$ denoted by
$A{\it Wr}B$ or $A\wr B$.

 Now, a presentation for a $\cal V$-verbal wreath product $A{\it Wr_V}B$ is
given based on the presentations of $A$ and $B$ as follows.

{\bf L}{\scr\bf EMMA} {\bf 2.7} [6]\ \ Let $A$ and $B$ be two groups with the following presentations:
$$ 1\lra R_1\lra F_1\st{\nu_1}{\lra}A\lra 1\ \ \ ,\ \ \ 1\lra R_2\lra
F_2\st{\nu_2}{\lra}B\lra 1\ . $$
Then    \\
 $$ 1\lra R\lra F\lra A{\it Wr_V}B\lra 1 $$
is a free presentation for the $\cal V$-verba wreath product
$A{\it Wr_V}B$, where

$(i)$ $F=F_1*F_2$ , the free product of $F_1$ and $F_2$ ;

$(ii)\ R=R_2S_V\ ;$

$(iii)\ S_V=[R_2,F_1]^FR_1^FR_V\ ;$

$(iv)\ R_V=V(H)\cap [(F_1^{f_2})^H]$ with $H=F_1^{F_2}$ (that is, $F_1^F$) and
$$ [(F_1^{f_2})^H]=<[x,y]|x\in (F_1^u)^H,y\in (F_1^v)^H,\ u,v\in F_2>\ .$$

 Now, the following theorem gives a structure of the Schur multiplier of a
$\cal V$-verbal wreath product.

{\bf T}{\scr\bf HEOREM} {\bf 2.8} [6]\ \ In accordance with the above notation, the following isomorphism holds.
 $$M(A{\it Wr_V}B)\cong M(A)\oplus M(B)\oplus \frac {C_V}{[C_V,
A{\it Wr}_*B]}\ \ .$$

{\bf C}{\scr\bf OROLLARY} {\bf 2.9}
$$M(A {\it Wr_*}B)\cong M(A)\oplus M(B)\ \ .$$
{\bf Proof.}\ \ We know that $A${\scr Wr}$_*B=A{\it Wr_V}B$,
where $V=\{1\}$. So $V(C)=1$ and then $C_V=1$. Now by Theorem 2.8
the result holds. $\Box $

{\bf C}{\scr\bf OROLLARY} {\bf 2.10}
 $$M(A\wr B)\cong M(A)\oplus M(B)\oplus \frac
 {[A_b]^*}{[[A_b]^*,A{\it
Wr_*}B]}\ \ ,$$
where $[A_b]^*$ is the cartesian subgroup of the free product $\prod_{b\in
B}^{*}A_b$ .\\
{\bf Proof.}\ \ As mentioned before $A\wr B=A{\it Wr_V}B$, where
$\cal V$ is the variety of abelian groups. Therefore $V(C)=C'$
and $C_V=V(C)\cap
[A_b]^*=[A_b]^*$ since $[A_b]^*\subseteq C'$. Now, the result holds by Theorem 2.8 . $\Box$\\

 We also need the following lemmas in our investigation.

{\bf L}{\scr\bf EMMA} {\bf 2.11} [16]\ \ Let $A$ be a subgroup of a group $G$, and $N$ be a normal subgroup of $G$ and
let $\{M_i|i\in I\}$ be a family of normal subgroups of $G$. Then
$$ [A\prod_{i}^{}M_i,N]=[A,N]\prod_{i}^{}[M_i,N]\ .$$

{\bf L}{\scr\bf EMMA} {\bf 2.12}

If $H$ is a subgroup of a finite abelian group $G$, then $G$ has
a subgroup isomorphic to $G/H$. Consequently, if A is a
homomorphic image of a finite
abelian group $G$, then $G$ has a subgroup isomorphic to $A$.\\

\begin{center}
{\bf 3. SOME RESULTS ON THE BAER INVARIANT OF A SEMIDIRECT
PRODUCT}
\end{center}

 First of all, by the functorial property of the Baer invariant, we are going to
generalize somehow Theorems 2.2 and 2.3 to an arbitrary variety.

{\bf T}{\scr\bf HEOREM} {\bf 3.1}\ \ Using the notion of the Baer
invariant, we can consider the following functor from the
category of all groups to the category of all abelian groups:
$${\cal V}M(-):{\cal G}{\it roups}\lra {\cal A}{\it b}\ ,$$
where $\cal V$ is an arbitrary variety of groups.

{\bf Proof.}\ \ See [5,9].

{\bf T}{\scr\bf HEOREM} {\bf 3.2}\ \ Let $G=B\st {\theta
}{\rhd\!\!\!<}A$ be the semidirect product of $A$ by $B$ under
$\theta $. Then ${\cal V}M(B)$ and ${\cal V}M(A)$ are direct
factors of ${\cal V}M(G)$, where $\cal V$ is an arbitrary variety.

{\bf Proof.}\ \ By the property of a semidirect product, we have the following
split exact sequence:
$$ 1\lra A\st {h}{\lra} G\st {f}{\lra }B\lra 1\ .$$
So there exists a homomorphism $g:B\ra G$ such that $f\circ
g=I_B$. Now, applying the functorial property of ${\cal V}M(-)$
(Theorem 3.1), the following exact sequence splits:
$$1\lra Ker({\cal V}M(f))\st {\subseteq }{\lra }{\cal V}M(G)\st {{\cal
V}M(f))}{\lra } {\cal V}M(B)\lra 1\ ,$$ and ${\cal V}M(f)\circ
{\cal V}M(g)=I_{{\cal V}M(B)}$. Hence, we have
$$ {\cal V}M(G)\cong {\cal V}M(B)\oplus Ker({\cal V}M(f))\ .$$
{\it i.e.} ${\cal V}M(B)$ is a direct factor of ${\cal V}M(G)$.

 To prove that ${\cal V}M(A)$ is a direct factor of ${\cal
 V}M(G)$, by a similar argument as above, we should consider the
 following split exact sequence:
$$1\lra {\cal V}M(A))\st {{\cal V}M(h) }{\lra }{\cal V}M(G)\st
{nat}{\lra } \frac{{\cal V}M(G)}{Im({\cal V}M(h))}\lra 1\ .\
\Box$$

 Now, in the following, we state a theorem of W. Haebich which is improves
somehow the structure of the complementary factor of $M(B)$ in $M(G)$ in Theorem 2.2 .

{\bf T}{\scr\bf HEOREM} {\bf 3.3} [6]\ \ By the  assumptions and notations of
Theorem 2.2, we have $$ M(G)\cong M(B)\oplus \frac {T\cap K'}{[T,K]}\ ,$$
where $K=F_1*B$ and
$$ T=<f_1^{-1}\ol{f_1}[b,f_1]\ |\ f_1,\ol{f_1}\in F_1;b\in
B;\nu_1\ol{f_1}=(\theta\nu_2b)(\nu_1(f_1))>^K\ .$$

 Now, we are going to find a structure similar to the above for the
complementary factor, ${\cal N}_cM(B)$ of ${\cal N}_cM(G)$ in Theorem 2.3 .

 First we need the following lemmas.

{\bf L}{\scr\bf EMMA} {\bf 3.4} [10]\ \ considering the assumptions and
notations of Lemma 2.1, the following statements hold:\\
$(i)\ R_1$ and $[R_2,F_1]$ are subgroups of $S$ ;\\
$(ii)\ R=R_2S$ ;\\
$(iii)\ R\cap \ga_{c+1}(F)=(R_2\cap \ga_{c+1}(F_2))(S\cap \ga_{c+1}(F))$, for
all $c\geq 1$ ;\\
$(iv)\ [R,\ _cF]=[R_2,\ _cF_2]\prod  [R_2,F_1,F_2]_c[S,\ _cF]$ , for
all $c\geq 1$ .

{\bf L}{\scr\bf EMMA} {\bf 3.5} [6]\ \ Let $C$ and $\ol C$ be two
normal subgroups of $A$ and let $A$ be a subgroup of a group $B$.
Let $\phi $ and $\ol \phi $ be two homomorphisms from $B$ to any
groups such that
$$ (A\cap Ker\phi )C=(A\cap Ker{\ol \phi}){\ol C}\ .$$
Then the map $\psi : \phi(A)/\phi (C)\lra {\ol \phi}(A)/{\ol \phi}({\ol C})$
given by
$$ (\phi a)\phi(C)\longmapsto ({\ol \phi}a){\ol \phi}({\ol C})$$
is an isomorphism.\\

 Now we are in a position to state and prove the main results of this section.

{\bf T}{\scr\bf HEOREM} {\bf 3.6}\ \ Let $G$ be a semidirect product of $A$ by $B$
under $\theta :B\ra Aut(A)$, and ${\cal N}_c$ be the variety of nilpotent
groups of class at most $c$. Then, using the notation of Theorem 3.3, there
exists an epimorphism as follows:
$$ {\cal N}_cM(G)\lra \!\succ {\cal N}_cM(B)\oplus \frac {T\cap \ga_{c+1}(K)}{
[T,\ _cK]}\ \ .$$
 So, in finite case, ${\cal N}_cM(G)$ has a subgroup isomorphic to
$$ {\cal N}_cM(B)\oplus \frac {T\cap \ga_{c+1}(K)}{[T,\ _cK]}\ \ .$$

{\bf Proof.}\ \ By the universal property of a free product, let $\delta $ be the natural
homomorphism from $F=F_1*F_2$ onto $K=F_1*B$ induced by $\nu_2:F_2\ra B$ and the
identity on $F_1$. If $b=\nu_2(f_2)$ and $(\theta b)(\nu_1(f_1))=\nu_1({\ol
f_1})$, then
$$\delta (f^{-1}_1{\ol f_1}[f_2,f_1])=\delta (f_1)^{-1}\delta ({\ol
f_1})[\delta (f_2),\delta (f_1)]=f^{-1}_1{\ol f_1}[b,f_1]\ .$$
Thus $\delta (S)=T$ and hence $\delta ([S,\ _cF])=[\delta(S),\ _c\delta(F)]=[T,\
_cK]$.\\
Also, we have
$$ \delta (S\cap \ga_{c+1}(F))\subseteq \delta(S)\cap \delta
(\ga_{c+1}(F))=T\cap \ga_{c+1}(K)\ .$$
To show the reverse containment, let $u\in T\cap \ga_{c+1}(K)$. Then there
exists $x\in S$ and $y\in \ga_{c+1}(F)$ such that $\delta(x)=u=\delta (y)$, and
so $yx^{-1}\in Ker\delta$. It is easy to see that $Ker\delta =R^F_2$, so
$yx^{-1}\in R^F_2$ and hence $y\in R^F_2S$.\\
By Lemma 3.4 parts $(i)$ and $(ii)$ we have $R=R^F_2S$. Therefore
$$ y\in R\cap \ga_{c+1}(F)\ .$$
Also, by Lemma 3.4 part $(iii)$ $R\cap \ga_{c+1}(F)=(R_2\cap
\ga_{c+1}(F_2))(S\cap \ga_{c+1}(F))$. Since $R_2\cap \ga_{c+1}(F_2)\subseteq
Ker\delta\ ,\ u=\delta (y)\in \delta (R\cap \ga_{c+1}(F))=\delta (S\cap
\ga_{c+1}(F))$, thus
$$ \delta (S\cap \ga_{c+1}(F))=T\cap \ga_{c+1}(K)\ .$$
Hence
$$ \frac {T\cap \ga_{c+1}(K)}{[T,\ _cK]}=\frac {\delta (S\cap
\ga_{c+1}(F)}{\delta ([S,\ _cF])}\ .$$
Now, putting $\phi =\delta :F\ra K$, ${\ol \phi}=identity:F\ra F,\ A=S\cap
\ga_{c+1}(F)$ and
$$ C=\prod [R_2, F_1, F_2]_c[S,\ _cF]\ \ ;$$
$$ {\ol C}=(\ga_{c+1}(F)\cap [R_2,F_1]^F)[S,\ _cF]\ \ ,$$
we have
$$ (A\cap Ker\phi )C=(S\cap \ga_{c+1}(F)\cap R^F_2)\prod [R_2,F_1,F_2]_c[S,\
_cF]$$
$$ =(S\cap \ga_{c+1}(F)\cap R^F_2)[S,\ _cF]\ \ \ \ \ \ \ (since\
\prod[R_2,F_1,F_2]_c\subseteq R^F_2\cap \ga_{c+1}(F)\ and\ S)$$
$$\ \ \ \ \ \ \ \ \ \ \ =(S\cap \ga_{c+1}(F)\cap R^F_2[R_2,F_1]^F)[S,\ _cF]$$
$$ =(S\cap \ga_{c+1}(F)\cap [R_2,F_1]^F)[S,\ _cF]\ \
(since\ \ F=F_2\rhd\!\!\!< F_1[F_1,F_2]\  and\ S\leq F_1[F_1,F_2])$$
$$ =(\ga_{c+1}(F)\cap [R_2,F_1]^F)[S,\ _cF]\ \ \ \ \ \ \ (By\
Lemma\ 3.4\ (i))$$
$$ (A\cap Ker{\ol \phi}){\ol C}\ .$$
So, by Lemma 3.5, we have the following isomorphism:
$$\frac {T\cap \ga_{c+1}(K)}{[T,\ _cK]}=\frac {\delta (S\cap
\ga_{c+1}(F))}{\delta ([S,\ _cF])}\cong \frac {S\cap
\ga_{c+1}(F)}{\ga_{c+1}(F)\cap [R_2,F_1]^F)[S,\ _cF]} $$
(Note that $\prod [R_2,F_1,F_2]_c\subseteq Ker\delta $). Clearly
$$\prod [R_2,F_1,F_2]_c[S,\ _cF]\subseteq (\ga_{c+1}(F)\cap [R_2,F_1]^F)[S,\
_cF]\ .$$
Therefore there exists a natural epimorphism
$$\frac {S\cap \ga_{c+1}(F)}{\prod [R_2,F_1,F_2]_c[S,\ _cF]}\lra \!\succ
\frac {S\cap \ga_{c+1}(F)}{(\ga_{c+1}(F)\cap [R_2,F_1]^F)[S,\ _cF]}\ .$$
 Now, by Theorem 2.3 there exists the following epimorphism
$$ {\cal N}_cM(G)\lra \!\succ {\cal N}_cM(B)\oplus \frac {T\cap \ga_{c+1}(K)}{
[T,\ _cK]}\ \ .$$
 Moreover, if $G$ is finite, then ${\cal N}_cM(G)$ is also
finite and hence by Lemma 2.12
$$  {\cal N}_cM(B)\oplus \frac {T\cap \ga_{c+1}(K)}{[T,\ _cK]}$$
is isomorphic to a subgroup of ${\cal N}_cM(G)$. $\Box$\\

 Now, using the above result, we can present an improved formula with respect to
that of Theorem 2.3  for the Baer invariant of a semidirect
product of an arbitrary group by a cyclic group with respect to
the variety ${\cal N}_c$.

{\bf T}{\scr\bf HEOREM} {\bf 3.7}\ \ Let $G$ be a semidirect product of an arbitrary group $A$ by a cyclic group $B$
under $\theta $ and ${\cal N}_c$ be the variety of nilpotent groups of class at
most $c$. Then, by the previous notation, there exists the following isomorphism
$$ {\cal N}_cM(G)\cong {\cal N}_cM(B)\oplus \frac {T\cap \ga_{c+1}(K)}{
[T,\ _cK]}\cong \frac {T\cap \ga_{c+1}(K)}{[T,\ _cK]} \ .$$

{\bf Proof.}\ \  Applying the previous proof, it is enough to show that
$$\prod [R_2,F_1,F_2]_c=\ga_{c+1}(F)\cap [R_2,F_1]^F\ .$$
Since $B$ is cyclic, we can consider the following presentation for $B$
$$1\lra R_2=<y^m>\lra F_2=<y>\st {\nu_2}{\lra }B\lra 1\ .$$
So $R_2$ has no commutators and hence by commutator calculus, specially on basic
commutators, we can conclude the required equality (see also [11] section 4).\\
Thus
$$ \frac {T\cap \ga_{c+1}(K)}{[T,\ _cK]} \cong
\frac {S\cap \ga_{c+1}(F)}{\prod [R_2,F_1,F_2]_c[S,\ _cF]}\ \ ,$$
and so, by Theorem 2.3, the result holds. $\Box$\\

{\bf N}{\scr\bf OTE} {\bf 3.8}\ \ In order to find a relation
between ${\cal N}_cM(A)$ and the
 quotient group $T\cap \ga_{c+1}(K)/[T,\ _cK]$, it is enough to
 consider the following natural epimorphism:
 $$ {\cal N}_cM(A)=\frac {R_1\cap \ga_{c+1}(F_1)}{[r_1,\ _cF_1]}\lra
 \!\succ \frac {(R_1\cap \ga_{c+1}(F_1))[T,\ _cK]}{[T,\ _cK]}\leq
 \frac {T\cap \ga_{c+1}(K)}{[T,\ _cK]} .$$
 Note that $R_1\leq T$ and $[R_1,\ _cF_1]\leq [T,\ _cK]$ and in
 finite case $T\cap \ga_{c+1}(K)/[T,\ _cK]$ has a subgroup which
 is isomorphic to a subgroup of ${\cal N}_cM(A)$.

  In the following theorem, we try to provide a considerable
  simplification of the quotient group $T\cap \ga_{c+1}(K)/[T,\
  _cK]$ and relate its structure more closely to that of $A$.

{\bf T}{\scr\bf HEOREM} {\bf 3.9}\ \ With the previous notation
and assumption the following isomorphism holds:
$$ \frac {U\cap \ga_{c+1}(D)}{[U,\ _cD]}\cong  \frac {T\cap
\ga_{c+1}(K)/ [T,\ _cK]}{(R_1\cap \ga_{c+1}(F_1))(\ga_{c+1}(K)\cap
[R_1,B]^K)[T,\ _cK]/[T,\ _cK]}\ \ , $$
 where $D=A*B$ and $U=<a^{-1}(\theta(b)(a))[b,a]|a\in A, b\in
 B>^D.$

{\bf Proof.}\ \ let $\eta : K=F_1*B\ra D=A*B$ be the natural
homomorphism induced by $\nu_1:F_1\ra A$ and the identity on $B$.
Clearly $Ker \eta =R_1^K$, $\eta (T)=U$ and hence $\eta ([T,\
_cK])=[U,\ _cD]$.\\
Also, similar to the proof of Theorem 3.6 we can show that
$$ \eta (T\cap \ga_{c+1}(K))= U\cap \ga_{c+1}(D)\ .$$

 Now, using Lemma 3.5 and a similar method of the proof of Theorem
 3.6, if we put $\phi =\eta :K\ra D$, ${\ol \phi}=identity:K\ra K,\ A=T\cap
\ga_{c+1}(K)$ and
$$ C=[T,\ _cK]\ \ ;$$
$$ {\ol C}=(R_1\cap \ga_{c+1}(F_1))(\ga_{c+1}(K)\cap [R_1,B]^K)[T,\ _cK]\ \ ,$$
then we have
$$ \frac {U\cap \ga_{c+1}(D)}{[U,\ _cD]}= \frac {\eta(T\cap \ga_{c+1}(K))}{\eta ([T,\ _cK])}
\cong  \frac {T\cap \ga_{c+1}(K)}{(R_1\cap
\ga_{c+1}(F_1))(\ga_{c+1}(K)\cap [R_1,B]^K)[T,\ _cK]}$$
$$ \cong
\frac {T\cap \ga_{c+1}(K)/ [T,\ _cK]}{(R_1\cap
\ga_{c+1}(F_1))(\ga_{c+1}(K)\cap [R_1,B]^K)[T,\ _cK]/[T,\ _cK]}.\
\  \Box $$

{\bf C}{\scr\bf OROLLARY} {\bf 3.10}\ \ If $A$ is a cyclic group,
then
$$ \frac {U\cap \ga_{c+1}(D)}{[U,\ _cD]}\cong  \frac {T\cap
\ga_{c+1}(K)}{[T,\ _cK]}\ \ .$$

 {\bf Proof.} Since $A$ is cyclic, we can consider $A=F_1/R_1$,
 where $F_1=<x>$ and $R_1=<x^n>$. So $R_1$ has no commutators and
 hence by commutator calculus we have
 $$ \ga_{c+1}(K)\cap [R_1,B]^K \subseteq [T,\ _cK]\ .$$
 Also, since $A$ is cyclic, so ${\cal N}_cM(A)=1$ and hence
 $$ R_1\cap \ga_{c+1}(F_1)=[R_1,\ _cF_1]\leq [T,\ _cK]\ .$$
 Thus, by Theorem 3.9 the result holds. $\Box$

{\bf C}{\scr\bf OROLLARY} {\bf 3.11}\ \ If $A$ and $B$ are cyclic
groups, then
$$ {\cal N}_cM(B\st {\theta
}{\rhd\!\!\!<} A)\cong \frac {U\cap \ga_{c+1}(D)}{[U,\ _cD]}\ .$$
{\bf Proof.} By Theorem 3.7 and Corollary 3.10 the result holds.
$\Box$

\begin{center}
{\bf 4. THE BAER INVARIANT OF A VERBAL WREATH PRODUCT}
\end{center}

 In this section, first, using a presentation for a verbal wreath product given
in [6], (Lemma 2.7), we are going to present a formula similar to
that of Theorem 2.3 for the Baer invariant of a ${\cal V}$-verbal
wreath product with respect to the variety of nilpotent groups,
${\cal N}_c$, where $\cal V$ is an arbitrary variety of groups.

 The following theorem is vital in the proof of the results of this section.

{\bf T}{\scr\bf HEOREM} {\bf 4.1}\ \ By the notation of Lemma 2.7, the following
equalities hold.

$(i)$ $R\cap \ga_{c+1}(F)=(R_2\cap \ga_{c+1}(F_2))(S_V\cap \ga_{c+1}(F))$, for
all $c\geq 1$;

$(ii)$ $[R,\ _cF]=[R_2,\ _cF_2]\prod [R_2,F_1,F_2]_c[S_V,\ _cF]$, for all $c\geq
1$,\\
where $\prod [R_2,F_1,F_2]_c$ was defined in Theorem 2.3 .

{\bf Proof.}\ \ First, by definition of $S_V$ and $R_V$ in Lemma 2.7, it is easy to
see that
$$ R_1\leq S_V\ ,\ [R_2,F_1]\leq S_V\ ,\ R_V\leq V(H)\leq
H=F_1^{F_2}=F_1[F_1,F_2]\ ,$$
$$S_V\leq F_1[F_1,F_2]\ ,\ F=F_2\rhd\!\!\!<F_1[F_1,F_2]\ \ and\ \ S_V\unlhd F\
\ .$$

 $(i)$ By Lemma 2.5
$$ \ga_{c+1}(F)=\ga_{c+1}(F_1)\ga_{c+1}(F_2)\prod [F_1,F_2]_{c+1}\ \ ,$$
where
$$ \prod [F_1,F_2]_{c+1}=<[F_1,F_2,F_{i_1},\ldots ,F_{i_{c-1}}]\ |\ i_j\in
\{1,2\},1\leq j\leq c-1>$$ and $\prod [F_1,F_2]_{c+1}\unlhd F$
(see also [13]). Now, by Lemma 2.4, we have
$$ R\cap \ga_{c+1}(F)=R_2S_V\cap
\ga_{c+1}(F_1)\ga_{c+1}(F_2)\prod [F_1,F_2]_{c+1} $$
$$ =(R_2\cap \ga_{c+1}(F_2))(S_V\cap
\ga_{c+1}(F_1)\prod [F_1,F_2]_{c+1}) $$
$$  =(R_2\cap \ga_{c+1}(F_2))(S_V\cap \ga_{c+1}(F))\ .  $$

$(ii)$ We use induction on $c$ . If $c=1$, then
$$ [R,F]=[R_2S_V,F]\ \ \ \  $$
$$\ \ \ \ \ \ \ \ \ \ \ \ \ \ =[R_2,F][S_V,F]\ \ \ , \ \ ({\it since\ }
S_V\unlhd F) $$
$$\ \ \ \ \ \ \ \ \ \ \ \ \ \ \ \subseteq [R_2,F_2][R_2,F_1]^F[S_V,F]\ \ , \
({\it since\ } F=<F_1,F_2>) $$
$$ \ \ \ \ \ \ \ \ \ \ \ \ \ \ \ \ \ \ =[R_2,F_2][R_2,F_1][S_V,F]\ \ \ , \ ({\it
since\ } [R_2,F_1]\leq S_V)\ \ . $$
 The reverse containment is clear. Hence $[R,F]=[R_2,F_2][R_2,F_1][S,F]$.

 Now, suppose $[R,\ _kF]=[R_2,\ _kF_2]\prod [R_2,F_1,F_2]_k[S_V,\ _kF]$.
Then we have
$$ [R,\ _{k+1}F]=[[R,\ _kF],F]$$    $$=[[R_2,\ _kF_2]\prod [R_2,F_1,F_2]_k[S_V,\
_kF],F]\ \ \ \ \ ({\it by\ induction\ hypothesis}) $$
$$ =[[R_2,\ _kF_2],F][\prod [R_2,F_1,F_2]_k,F][[S_V,\ _kF],F]$$
$$({\it since}\ [S_V,\ _kF]\ ,\ \prod[R_2,F_1,F_2]_k\unlhd F)  $$
$$ \subseteq [[R_2,\ _kF_2],F_2]\prod [R_2,F_1,F_2]_{k+1}[S_V,\
_{k+1}F] $$
$$ \subseteq [R,\ _{k+1}F]\ .$$
Therefore, by induction we have
$$ [R,\ _cF]=[R_2,\ _cF_2]\prod [R_2,F_1,F_2]_c[S_V,\ _cF]\ \ \ \ for\ all\
c\geq 1\ .\ \ \Box $$

 Now, we are in a position to prove one of the main result of this section.

{\bf T}{\scr\bf HEOREM} {\bf 4.2}\ \ Let $G$ be the $\cal V$-verbal wreath product of $A$ by $B$, where $\cal V$ is
an arbitrary variety of groups. Then, using the previous notation, the following
isomorphism exists.
$$ {\cal N}_cM(G)={\cal N}_cM(A{\it Wr_V}B)\cong {\cal N}_cM(B)\oplus \frac {
S_V\cap \ga_{c+1}(F)}{\prod [R_2,F_1,F_2]_c[S_V,\ _cF]}\ .$$

{\bf Proof.}\ \ By the previous assumptions and notations the following natural homomorphisms
exists
$$F\st{\varphi}{\lra} \frac{F}{[R_2,\ _cF_2]^F}\st{\eta}{\lra}\frac{F}{[R_2,\
_cF_2]\prod [R_2,F_1,F_2]_c[S_V,\ _cF]}\ \ . $$
 By Theorem 4.1 the following isomorphisms hold:
$$ \frac{R\cap \ga_{c+1}(F)}{[R,\ _cF]}\cong (\eta \varphi )(R\cap
\ga_{c+1}(F)) $$
$$ \cong (\eta \varphi )(R_2\cap \ga_{c+1}(F_2))(\eta
\varphi )(S_V\cap \ga_{c+1}(F))\ .\ \ \ \ \ (*)  $$ Consider the
following two natural homomorphisms:
$$\frac{F_1*F_2}{[R_2,\ _cF_2]^F} \st{h}{\lra}  F_1*\frac{F_2}{[R_2,\ _cF_2]}
\st{g}{\lra} \frac{F_1*F_2}{[R_2,\ _cF_2]^F}\ \ \ \ \ , $$
given by
$$\ol{f_1}  \longmapsto   f_1\ \ \ \ \ \ \ \ \ \ f_1  \longmapsto  \ol{f_1}
\ \ \ \ $$
$$ \ol{f_2} \longmapsto  \ol{f_2}\ \ \ \ \ \ \ \ \ \ \ol{f_2} \longmapsto
\ol{f_2} \ \ \ . $$
It is easy to see that $h\circ g=1$ and $ g\circ h=1 $ i.e $h$ and $g$ are
isomorphisms. Thus, we have
$$ \frac{F_1*F_2}{[R_2,\ _cF_2]^F}=\varphi (F)\cong F_1*\frac{F_2}{[R_2,\
_cF_2]}\ \ \ .$$
Also
$$\varphi (F_2)=\frac {F_2[R_2,\ _cF_2]^F}{[R_2,\ _cF_2]}=\frac {F_2[R_2,\
_cF_2][R_2,\ _cF_2,F_1]^F}{[R_2,\ _cF_2][R_2,\ _cF_2,F_1]^F}$$
$$\cong \frac {F_2}{F_2\cap [R_2,\ _cF_2][R_2,\ _cF_2,F_1]^F}\ \ \ (By\ the\
second\ isomorphism\ Theorem)$$
$$=\frac {F_2}{[R_2,\ _cF_2]}\ \ \ \ \ \ \ \ \ (By\ Lemma\ 2.5)$$
and
$$\varphi (F_1[F_1,F_2])\cong
\varphi (F_1)[\varphi (F_1),\varphi (F_2)]\cong F_1[F_1,\frac {F_2}{[R_2,\
_cF_2]}]\ \ . $$
(Note that
$$ \varphi (F_1)=\frac {F_1[R_2,\ _cF_2]^F}{[R_2,\ _cF_2]^F}\cong \frac
{F_1}{F_1\cap [R_2,\ _cF_2]^F}=\frac {F_1}{1}\cong F_1\ .)$$
Therefore
$$ \varphi (F)\cong F_1*\frac{F_2}{[R_2,\ _cF_2]}\cong \frac{F_2}{[R_2,\
_cF_2]}\rhd\!\!\! <F_1[F_1,\frac{F_2}{[R_2,\ _cF_2]}] $$
 $$ \ \ \ \ \ \ \ \ \ \ \ \ \ \ \ \ \ \ \ \ \ \ \ \ \ \
\cong \varphi(F_2)\rhd\!\!\!<\varphi(F_1)[\varphi(F_1),\varphi(F_2)]\ \ .$$
Clearly $Ker(\eta )=\varphi (\prod [R_2,F_1,F_2]_c[S_V,\ _cF])$
$$\leq \varphi (F_1)[\varphi (F_1),\varphi (F_2)]\ . \ \ (**)$$
So by [6, (Lemma 1.3)] we have
$$ (\eta \varphi )(F)\cong \frac{\varphi(F)}{Ker(\eta )}=\frac{\varphi(F)}{
\varphi(\prod [R_2,F_1,F_2]_c[S_V,\ _cF])}$$
$$\cong \varphi(F_2)\rhd\!\!\! <\frac{\varphi(F_1)[\varphi(F_1),\varphi(F_2)]}
{Ker(\eta )}\ \ . $$
By $(**)$, it is easy to see that  $(\eta \varphi )(F_2)\cong \varphi(F_2)$ and
$(\eta\varphi )(F_1)\cong \varphi(F_1)/Ker(\eta  )$.
Thus we have
$$ (\eta \varphi )(F)\cong
\varphi(F_2)\rhd\!\!\! <\frac{\varphi(F_1)[\varphi(F_1),\varphi(F_2)]}{Ker(\eta
)} $$
$$\ \ \ \ \ \ \ \ \ \ \ \cong (\eta\varphi )(F_2)\rhd
\!\!\! <(\eta\varphi )(F_1)[(\eta\varphi )(F_1),(\eta\varphi
)(F_2)]\ \ \ .$$ So, we can conclude that
$$(\eta \varphi )(R_2\cap \ga_{c+1}(F_2))\ \cap\ (\eta\varphi )(S_V\cap
\ga_{c+1}(F))$$      $$  \subseteq (\eta\varphi )(F_2)\cap
(\eta\varphi )(F_1)[(\eta\varphi )(F_1),(\eta\varphi )(F_2)]=1 $$
Hence, by considering $(*)$, we have
$$ \frac{R\cap \ga_{c+1}(F)}{[R,\ _cF]}\cong (\eta\varphi )(R_2\cap
\ga_{c+1}(F_2))\oplus(\eta\varphi )(S_V\cap \ga_{c+1}(F))\ . $$
On the other hand, we have the following isomorphisms.
$$ (\eta \varphi )(R_2\cap\ga_{c+1}(F_2))=\frac{(R_2\cap
\ga_{c+1}(F_2))Ker(\eta \varphi )}{Ker(\eta \varphi )}$$
$$\cong \frac{R_2\cap \ga_{c+1}(F_2)}{(R_2\cap \ga_{c+1}(F_2))\cap Ker(\eta
\varphi )}\ \ \ (by\ the\ second\ isomorphism\ Theorem) $$
$$\ \ \ \ \ \ \ \ \cong \frac{R_2\cap \ga_{c+1}(F_2)}{[R_2,\ _cF_2]}\ \ \ (by\
Lemma\ 2.5)$$
$$\cong {\cal N}_cM(B)\ \ .$$
Also
$$ (\eta \varphi )(S_V\cap \ga_{c+1}(F))=\frac{(S_V\cap
\ga_{c+1}(F))Ker(\eta\varphi )}{Ker(\eta\varphi )}$$
$$\cong \frac{S_V\cap \ga_{c+1}(F)}{(S_V\cap \ga_{c+1}(F))\cap Ker(\eta\varphi
)} \ \ \ \ (by\ the\ second\ isomorphism\ Theorem)$$
$$\ \ \ \ \ \ =\frac{S_V\cap \ga_{c+1}(F)}{\prod [R_2,F_1,F_2][S_V,\ _cF]}\ \ \
\ (by\ Lemma\ 2.5)\ .$$
Therefore
$$ {\cal N}_cM(B{\it Wr_V}A)\cong \frac{R\cap \ga_{c+1}(F)}{[R,\ _cF]}$$ $$\cong
{\cal N}_cM(B)\oplus \frac{S_V\cap \ga_{c+1}(F)}{\prod [R_2,F_1,F_2]_c[S_V,\
_cF]}\ \ .\ \ \ \Box $$

 Now, we can obtain the following corollary if we put $c=1$.

{\bf C}{\scr\bf OROLLARY} {\bf 4.3} [6]\ \ Suppose $G=A{\it
Wr_V}B$ is the $\cal V$-verbal wreath product of $A$ by $B$,
where $\cal V$ is a variety of groups. Then
$$ M(A{\it Wr_V}B)\cong M(B)\oplus \frac {S_V\cap F'}{[R_2,F_1][S_V,F]}\ . $$

 Now, in what follows, we are going to improve the structure of the
complementary factor of ${\cal N}_cM(B)$ in ${\cal N}_cM(A{\it
Wr_V}B)$.

{\bf N}{\scr\bf OTATION} {\bf 4.4}\ \
 By considering the previous notation, let $\delta :F_1*F_2\lra F_1*B$ be the
natural epimorphism, where $A\cong F_1/R_1$ and $B\cong F_2/R_2$
are presentations for $A$ and $B$. Put $T_V=\delta (S_V)$ and
$K=F_1*B$. Using Lemmas 3.4 and 2.7, Theorem 4.1 and similar
properties of $S_V$ to that of $S$, we can rewrite the proof of
Theorem 3.6 to get the following Theorem about the Baer invariant
of a verbal wreath product.

{\bf T}{\scr\bf HEOREM} {\bf 4.5}\ \ By the notation above, there
exists the following epimorphism:
$$ {\cal N}_cM(A {\it Wr_V}B)\lra \!\succ {\cal N}_cM(B)\oplus \frac { T_V\cap
\ga_{c+1}(K)}{[T_V\ _c,K]}\ .$$ Also, in finite case, ${\cal
N}_cM(A{\it Wr_V}B)$ has a subgroup isomorphic to
$$ {\cal N}_cM(B)\oplus \frac { T_V\cap \ga_{c+1}(K)}{[T_V\ _c,K]}\ .$$

 Now, by putting a condition on $B$, and a similar proof to that of Theorem 3.7,
we can present a formula for a verbal wreath product which is somehow better
than that of Theorem 4.2 .

{\bf T}{\scr\bf HEOREM} {\bf 4.6}\ \ By the above notation, let
$G$ be a $\cal V$-verbal wreath product of an arbitrary group $A$
by a cyclic group $B$, where $\cal V$ is any variety of groups.
Then
$$ {\cal N}_cM(A {\it Wr_V}B)\cong
 {\cal N}_cM(B)\oplus \frac { T_V\cap \ga_{c+1}(K)}{[T_V\ _c,K]}
= \frac { T_V\cap \ga_{c+1}(K)}{[T_V\ _c,K]}\ .$$

 {\bf N}{\scr\bf OTE} {\bf 4.7}\ \ Note that W. Haebich in [6] prove that
$$ M(A{\it Wr_V}B)\cong M(B)\oplus \frac { T_V\cap K'}{[T_V,K]}\ \ ,$$
where $A$ and $B$ are two arbitrary groups.\\

\begin{center}
{\bf 5. THE BAER INVARIANT OF A FREE WREATH PRODUCT}
\end{center}

 In this final section we try to find a structure for the Baer invariant of a
free wreath product with respect to the variety of nilpotent
groups. We recall that a free wreath product is in fact a $\cal
V$-verbal wreath product where $\cal V$ is the variety of all
groups, i.e. $V=\{1\}$. Then, by the previous notation, Lemma
2.7, we have $V(H)=1$ and so $R_V=1$ and $S_V=[R_2,F_1]^FR_1^F$.
Now, by the definition of the natural epimorphism $\delta :
F_1*F_2\lra F_1*B$, $\delta (R_2)=1$ and hence $T_V=\delta
(S_V)=R_1^K=R_1[R_1,K]^K$.

{\bf L}{\scr\bf EMMA} {\bf 5.1}\ \ By the above notation, we have

$(i)$  $[T_V,\ _cK]=[R_1,\ _cF_1]D_c$,\\ where $D_c=\prod_{\st
{G_i\in \{F_1,B\}}{ \exists j\ G_j=B}}[R_1,G_1,\ldots ,G_c]^K$;

$(ii)$ $T_V\cap \ga_{c+1}(K)=(R_1\cap \ga_{c+1}(F_1))E_c$,\\
where $E_c=[R_1,K]^K\cap \ga_{c+1}(K)\cap [F_1,B]$.

{\bf Proof.}

$(i)$ $[T_V,\ _cK]=[R_1,\ _cK]^K=[R_1,\ _K]^K=[R_1,\ _cF_1]D_c$\
\  $(by\ commutator\\ calculus\ and\ Lemma\  2.11)$.

$(ii)$ $T_V\cap \ga_{c+1}(K)=R_1^K\cap
\ga_{c+1}(K)=R_1[R_1,K]^K\cap \ga_{c+1}(K)$
\\ $\subseteq (R_1\cap \ga_{c+1}(F_1))([R_1,K]^K\cap \ga_{c+1}(K))$ \ \ $(by\
Lemmas\ 2.4\ and\ 2.5)$
$$ =(R_1\cap \ga_{c+1}(F_1))E_c\ \ .$$
The reverse inclusion can be seen easily. $\Box $

{\bf L}{\scr\bf EMMA} {\bf 5.2}\ \ By the above notation, we have
the following isomorphism:
$$\frac {T_V\cap \ga_{c+1}(K)}{[T_V,\ _cK]}\cong {\cal N}_cM(A)\oplus \frac {E_c}{D_c}\ \ \ .$$

{\bf Proof.}\ \ Put
$$ \varphi : K \st {nat}{\lra \!\succ} \frac {K}{[T_V,\ _cK]} \ \ \ and $$
$$K=F_1*B\st{\xi=nat}{\lra\!\succ} \frac {F_1}{[R_1,\
_cF_1]}*B\st{\rho}{\lra \! \succ} \frac {\xi(K)}{\xi(D_c)}\ . $$
Clearly $ker(\rho\xi)=[R_1,\ _cF_1]D_c=[T_V,\ _cK]=ker\varphi $,
so $\varphi=\rho\xi $. Now similar to the proof of Theorem 3.9 in
[6] we can show that
$$ \varphi(T_V\cap \ga_{c+1}(K))=\varphi(R_1\cap
\ga_{c+1}(F_1))\oplus \varphi(E_c) $$ and
$$\varphi(R_1\cap \ga_{c+1}(F_1))\cong \xi (R_1\cap\ga_{c+1}(F_1))\cong
{\cal N}_cM(A)\ \ .$$ Also
$$\varphi(E_c)=\frac {E_c([R_1,\ _cF_1]D_c)}{[R_1,\
_cF_1]D_c}\cong \frac {E_c}{[R_1,\ _cF_1]D_c\cap E_c}=\frac
{E_c}{([R_1,\ _cF_1]\cap E_c)D_c}\ \ .$$ Since $[R_1,\ _cF_1]\cap
E_c=1$, so $\varphi(E_c)\cong E_c/D_c$ . Hence the result holds.
$\Box $

 Now, we are in a position to find a structure for the $c$-nilpotent multiplier
of a free wreath product.

{\bf T}{\scr\bf HEOREM} {\bf 5.3}\ \ Let $G$ be the free wreath
product of $A$ by $B$ i.e. $G=A{\it Wr_*}B$. Then there exists an
epimorphism as follows:
$${\cal N}_cM(A{\it Wr_*}B)\lra \! \succ {\cal N}_cM(A)\oplus {\cal N}_cM(B)
\oplus \frac {E_c}{D_c}\ \ .$$ Also, in finite case, the above
structure is isomorphic to a subgroup of ${\cal N}_cM(A{\it
Wr_*}B)$.

{\bf Proof.}\ \ By Theorem 4.5 and Lemma 5.2 the result holds.
$\Box $

{\bf T}{\scr\bf HEOREM} {\bf 5.4}\ \ Suppose $A$ and $B$ are two
groups such that $B$ is cyclic. Then
$${\cal N}_cM(A{\it Wr_*}B)\cong {\cal N}_cM(A)\oplus {\cal N}_cM(B)\oplus \frac {E_c}{D_c}\cong
{\cal N}_cM(A)\oplus \frac {E_c}{D_c}\ \ .$$

{\bf Proof.}\ \ By Theorem 4.6, Lemma 5.3 the result holds. $\Box$

\end{document}